\input amstex

\documentstyle{amsppt}

\magnification = \magstep 1
\hsize=6.5 true in
\vsize=8.9 true in

\define\rank{\operatorname{rank}}

\define\bho{H^0(C,}
\define\ho{h^0(C,}

\define\hi{h^1(C,}
\define\rnk{\operatorname{rank} } 
\define\mua{\mu_\alpha}
\define\red{\operatorname{red} }

\topmatter

\title 
Birational maps of moduli of  Brill-Noether pairs
\endtitle

\author David C. Butler \endauthor


\address School of Mathematics, 
Tata Institute of Fundamental Research, 
Homi Bhabha road, Bombay, India  \endaddress

\address e-mail: butler\@math.tifr.res.in \endaddress

\abstract
Let $C$ be a smooth projective irreducible curve of genus $g$.
And let $G_{\alpha}(n,d,l)$ be the moduli space of $\alpha$
stable pairs of a vector bundle of $\rank n$, $\deg d$ and a subspace
of $\bho E)$ of $\dim = l $.  We find an explicit birational
map from $G_{\alpha} (n, d, n+1)$ to $G_{\alpha} (1, d, n+1)$
for $C$ general, $\frac{1}{\alpha} \gg 0$ and $g \ge n^2-1$.
Because of this and other examples, we conjecture 
$G_{\alpha} (a, d, a+z)$ maps birationally to
$G_{\alpha} (z, d, a+z)$ for $\frac{1}{\alpha} \gg 0$
and $C$ general with $g>2$.
\endabstract

\endtopmatter

\document

\baselineskip = 20 pt
\lineskip = 11 pt
\lineskiplimit = 10 pt


\heading Introduction \endheading

Let $C$ be a smooth projective irreducible curve of genus $g$ 
defined over $\Bbb C$.
And for a vector bundle $E$ over $C$ let $\mu(E)=\frac{\deg (E)}
{\rank (E)}$.  We say $E$ is stable (or semistable) if for
every subbundle $S\subsetneq E$ we have $\mu(S)< \mu(E)$
(or $\le$).  The chief advantage of stable bundles is that
they form a coarse moduli space $M=M^{(s)}(n,d)$ of 
$\rank n$, $\deg d$ stable bundles.  Contained in M is the
Brill-Noether locus $W=W^r_{n,d}$ of $E\in M$ such that
$\ho E)\ge r+1$.  It is also known that $W$ is locally 
a determinantal variety with `theoretical' dimension given by 
$\rho (r, n, d) = n^2(g-1) +1 -(r+1)(r+1 -d +n(g-1) )$.
And hence, if $W$ is not empty it has the above minimal dimension.
It is easily seen that $\rho (a+z-1, a, d) =\rho (a+z-1, z, d)$.
This suggest some sort of relationship between
$W^{a+z-1}_{a,d}$ and $W^{a+z-1}_{z,d}$ --- perhaps
a birational isomorphism.  A candidate for that birational
isomorphism is the following.  Let $E\in W^{a+z-1}_{a,d}$ be generated by
global sections, and consider the sequence of vector bundles 
with $M_E$ the kernel of the evaluation map.
$$0 \to M_E \to H^0(C,E)\otimes \Cal O_C \to E \to 0.$$
Now if $E$ has no trivial summands $H^0(C,E)^*\subseteq \bho M_E^*)$.
So we hope $M_E^*\in W^{a+z-1}_{z,d}$.  Furthermore, if
$H^0(C,E)^* =  \bho M_E^*)$, then there is an inverse map
using the dual of the exact sequence above.  
This is a special case of the
dual span map which we will define shortly.

There are problems with the map $E$ goes to $M_E^*$.
\roster 
   \item $E$ may not be spanned.
   \item $M_E$ may not be stable.
   \item $\ho E)$ may be $> a+z$.
   \item $\ho M_E^*)$ may be $> a+z$.
\endroster
Problems 1 and 2 are serious.  They are solved under some conditions
in this paper and are conjectured here to be solvable under `most' 
conditions.  Sometimes they fail and nothing can be done.  
But 3 and 4 have been solved by Raghavendra and Vishwanath.  
Rather than looking at bundles, we
look at Brill-Noether pairs $(E,V)$ which are pairs of a vector bundle
$E$ and a space of sections $V\subseteq \bho E)$.  
A pair $(E,V)$ is of type $(n,d,l)$ if $\rank(E)=n$, $\deg(E)=d$
and $\dim (V) = l$.

A morphism of Brill-Noether pairs 
$(F,W) \to (E,V)$ is a morphism $F\to E$ such that the natural
morphism $W \to \bho E)$ factors through $V$.  A morphism is a
subbundle if $F$ is a subbundle of $E$.  One should note
that if $F = E$, an automorphism of the bundle may induce an
isomorphism of pairs $(E, V)$ and $(E, W)$ with $V\ne W$
as subspaces of $\bho E)$.  This never happens to simple
bundles (such as stable bundles) because the only automorphisms
are scalars.

To form a moduli space of pairs requires a notion of stability.
Following King and Newstead \cite{3}, we choose a rational
number $\alpha > 0$.  Then we define slope.
$$\mu_\alpha = \frac{\deg (E) + \alpha \dim(V)}{\rnk(E)}.$$
Now we have the usual definition of stability.  $(E,V)$ is
$\alpha$ stable (or $\alpha$ semistable) if for every subbundle $(F,W)$,
$\mua (F,W) < \mua (E,V)$ (or $\le$).

\subheading{Definition}  The set of Brill-Noether
pairs of type $(n,d,l)$ with $\alpha$ stability
is $G_{\alpha} (n,d,l)=$ the set of isomorphism
classes of $\alpha$ stable pairs $(E,V)$ with $\rank(E) = n$,
$\deg(E) =d$ and $\dim (V)=l$.
And if no $\alpha$ appears we are considering the set
(of not necessarily stable) Brill-Noether pairs.

We will use $\alpha$ with $\frac{1}{\alpha} \gg 0$.  To see how this
works assume $\rank E =n$ is fixed and $\ho E)=h$ is fixed.  
For $F$ a subbundle of $E$, $\mu(E)-\mu (F)\ge \frac{1}{n^2}$.  
So set $\alpha < \frac{1}{n^2}/ h$.
Now if $E$ is stable, $(E,V)$ is stable.  And if $E$ is unstable
the pair $(E,V)$ is unstable.  But if $E$ is semistable, three things can
happen.  Suppose $\mu(F) =\mu(E) $ and $F\subsetneq E$.  and $(F,W)$ is a
subpair of$(E,V)$.  If for all such $F$, 
$\frac{\dim W}{\rank F} < \frac{\dim V}{\rank E}$ (or $\le$), 
then the pair $(F,W)$ is $\alpha$ stable (or $\alpha$ semistable).  
Otherwise the pair is $\alpha$ unstable. 

Raghavendra and Vishwanath \cite {6} construct a moduli
space for Brill Noether pairs with $\alpha$ small.
And even though we do not use it, we note
King and Newstead \cite {3} have constructed a coarse moduli space
of $\alpha$ stable Brill Noether pairs for $\alpha> 0$ and
rational on a polarised curve (which need not be reduced).

Now back to the dual span.  
\subheading{Definition}
The set of spanning pairs is the set of pairs $(E,V)$,
where $V$ spans $E$ and has no trivial summands.
The set of $\alpha$ stable spanning pairs is:
$$S_\alpha(n,d,l) =\{(E,V)\in G_\alpha (n,d,l)
\text{ such that }  E \text{ is spanned by }V\}$$.

Now we revisit the dual span map.  Let $(E,V)$ be a spanning pair
which contains no trivial summand (but is not necessarily stable).
We have an exact sequence:
$$0\to M_{V,E} \to V\otimes \Cal O_C \to E \to 0.$$
where $M_{V,E}$ is just the kernel of the evaluation map.  
And the dual gives us:
$$0 \to E^* \to V^*\otimes \Cal O_C \to M_{V,E}^* \to 0.$$
Notice the pair $(M_{V,E}^*,V^*)$ is a (not necessarily
stable) spanning pair.  So this map needs a name.

\subheading{Definition}
The dual span map is the bijection of the set of
spanning pairs which takes $(E,V)$ to $(M_{E,V}^*, V^*)$,
the dual of the kernel of the evaluation map.
If $V=\bho E)$ we may write $M_E$.

If the type of $(E,V)$ is $(n,d,l)$ then the type of the 
dual span is $(l-n,d,l)$.
Furthermore, the dual span of $(M_{V,E}^*, V^*)$ is just $(E,V)$.
The problem then is stability.  If $(V,E)$ is $\alpha$ stable,
is its dual span $\alpha$ stable?  Not always.  But perhaps often.
And is a general $\alpha$ stable Brill-Noether pair a spanning pair 
for a given type $(n,d,l)$ with $l>n$?
Where by general we mean an Zariski open set which is
dense in each irreducible component.

Now the main result.
\proclaim {Main Theorem}
If $C$ is a general curve, $\frac{1}{\alpha} \gg 0$, and
$g\ge n^2-1$ or $g=n=2$, then
$$
S_{\alpha}(1,d,n+1)_{\red} \cong S_{\alpha} (n,d,n+1)_{\red}.
$$
\endproclaim

\subheading{Remark 1}  
We could state a stronger result for $n=2$ because Teixidor i Bigas \cite{8}
has proven that $W^3_{2,d}$ is reduced and irreducible.
And Tan \cite{7} has proven the locus is non-empty.
The problem with extending this to the case
$n>2$ is showing a general pair in
the space $G_{\alpha} (n,d,n+1)$ is spanned.  Although we show
some component is generally spanned.  
Nothing is known about the scheme structure in general.
And the condition $g\ge n^2-1$ seems to be a flaw in our
proof and not a part of nature.

\proclaim{Conjecture 1}
Fix $(n,d,l)$ with $l>n$, let $C$ be a general curve with $g>2$,
and choose $\alpha$ so $\frac{1}{\alpha} \gg 0$.
Then $S_{\alpha} (n,d,l)$ is dense in $G_{\alpha} (n,d,l)$.
\endproclaim

\proclaim{Conjecture 2}
Fix $(n,d,l)$ with $l>n$ and let $C$ be a general curve with $g>2$.
For a general $(E,V) \in S_{\alpha} (n,d,l)$ 
with $\frac{1}{\alpha}\gg 0$, $M_{V,E}$ is $\alpha$ stable.
and the morphism $S_{\alpha} (n,d,l) \to S_{\alpha} (l-n,d,l)$
is birational (at least after reducing the schemes).
\endproclaim

\heading \S 1 Stability Results \endheading

Now we prove that in some cases, the bundle $M_E$
is stable.  A trick for $\rank 1$ allows us to show
some $M_{L,V}$ with $V\ne \bho L)$ are stable on a general curve.  
We shall also indicate why we need general bundles on general curves.

The first result that needs mentioning is that if $E$
is stable and $\mu(E)> 2g$ then $M_E$ is stable \cite {1, theorem 1.2}.  
If $l>gn$ then we get a morphism:
$$G_{\alpha}(n, l+gn, l+n) \to G_{\alpha}(l,l+gn, l+n).$$
In \cite {5}, Mercat proves the above is an isomorphism.
This beautiful result provides evidence for the conjecture.
But this should not mislead the reader to believe we
will have an isomorphism in general.  Or that we should
get a birational map on special curves.

It can be shown 
\cite {1, proposition 1.5 and example 2.6}
that for a small number $\epsilon >0$ there is a bundle
(on any curve) with $\mu(E)\ge 2g - \epsilon$ and $M_E$ unstable.
This necessitates the conjecture's assumption that we deal with
a general pair $(E,V)$.

It is also assumed that the curve is general.  Suppose
to the contrary that $C$ is special, in fact
hyperelliptic, with hyperelliptic
bundle $A$.  If $\mu(E)>g+1$ then $\ho A^*\otimes E)>0$.
So $A \subset E$ and $M_A \subset M_E$. Now
$\mu(M_A)=-2$, and if $\hi E)=0$ and $\mu(E) < 2g$,
we have $\mu(M_E)=\frac{-\mu(E)}{\mu(E)-g} <-2$.
So $M_E$ is unstable.  A similar construction applies
to any fixed gonality $\beta$ when $g \gg 0$.

Now we get a more systematic theorem for $E=L$ a line
bundle on a general curve.

\proclaim{Theorem 2}
Let $C$ be a smooth projective irreducible curve of genus $g$.
And assume that if $\rho (r,d,1) < 0$, then $W^r_{d,1}= \emptyset$.
If $L$ is generated by sections $M_L$ is semistable.  In fact,
it fails to be stable iff {\bf all\/} the following hold.
\roster
\item $\hi L) = 0$.
\item $\deg L = g + r$ and $r|g$.
\item There is an effective divisor $Z$ with 
      $\ho L(-Z)) = \ho L)-1$ and $\deg Z = 1+\frac{g}{r}$.
\endroster
\endproclaim

\subheading{Remark 2}
Assuming that C has the property that $\dim (W^r_d)=\rho(g,r,d)$.
If $g\ge 3$ then by a dimension count, a general $L$ satisfying 1 and 2
does not satisfy 3.  And for any genus a general $L\in W^r_d$ but
not in $W^{r+1}_d$ is spanned by a dimension count.

\demo{Proof of Theorem 2}
Given $L$ generated by sections and $S$ a subbundle of $M_L$,
there is a commutative diagram:
$$
\CD
@.     0       @.          O          @.              \\
@.   @VVV                 @VVV                        \\
0 @>>> S @>>> V\otimes \Cal O_C @>>> E @>>> 0         \\
@.   @VVV                 @VVV      @VVV              \\
0 @>>> M_L @>>> \bho L)\otimes \Cal O_C @>>> L @>>> 0 \\
\endCD
$$
$S= M_{E,V}$.  It would be simpler if $E$ was a line bundle.
It is possible to take hyperplane sections and get $S=M_{W,A}$
where $A$ is the determinant of $E$.  But there is an obvious
question.  Is $W=\bho A)$?  The answer is generally no.
If $E=A$ (and hence $A$ is a subbundle of $L$), then the answer
is yes.  But suppose $E\ne A$.  In taking hyperplane sections,
we get 
$$0\to \oplus \Cal O_C \to E \to A \to 0.$$
And this sequence is exact on global sections.  Now if we
tensor by $\omega_C$ we get a sequence which is not exact
on global sequences.  This is because $E$ has no trivial
summands and hence no degree $0$ quotient bundles, and hence
$E\otimes \omega_C$ has not quotient bundles of degree $\le 2g-2$.
That means $E\otimes \omega_C$ is non special.  But the subbundles
$\omega_C$ are special.

The failure of the seqeunce tensored by $\omega_C$ means
$\bho A) \otimes \bho \omega_C) \to \bho A\otimes \omega_C)$
does not surject.  But by a theorem of Green \cite{2, Theorem 4.b.2}
this only happens if the image of the morphism induced by $A$
is a rational normal curve.  If it is a pencil the base point
pencil trick applies.  So it must have image of rational normal
curve of $\deg \ge 2$.  Given our assumption on generality of
the curve, the gonality of the curve is $\ge \frac{g+2}{2}$
so the only possibility is if $\deg (A) = g+2$ and $\ho A)=3$
and hence $A$ is nonspecial.  So a subbundle of the form $M_A$
for $A$ a special bundle must be induced by a subbundle $A\subseteq L$.
And all other special bundles are of the form $M_{A,W}$ for $W$ not complete.

Now suppose $\ho L)=r+1$.  Then a simple calculation shows
$\mu(M_L)\ge \frac{-g-r}{r}$ with equality iff $L$ is nonspecial.
Now if $M_L$ is destabilized by a subbundle of the form $M_{A,W}$,
$\dim (W) = s \le r-1 $ and $\deg(A) < \deg (L)$.  
Since $\mu(M_{A,W})= -\deg (A)/(\dim (W) -1)$,
we want $\deg (A) = \delta$ as small as possible for a given $s$.  
So we use the Brill-Noether numbers.
We have:
$$
\align
g-(s+1)(s+1 -\delta +g-1) &\ge 0 \\
\delta &\ge s + g -\frac{g}{s+1}.
\endalign
$$
If the last inequality is exact, then $(s+1)|g$.
Now we claim 
$$
\align
\mu(M_L) \ge \frac{-g-r}{-r} &\ge \frac{-s -g}{s} + \frac{g}{(s+1)s} \ge 
                                                       \mu(M_{A,V})   \\
            \frac{-g}{r} &\ge \frac{-g}{s} + \frac{g}{s(s+1)}         \\
          -g(s(s+1)) &\ge -gr(s+1) +gr                                \\
          -gs^2 -gs  &\ge -gs^2 -g(r-s)s.
\endalign
$$
So $\mu(M_L) \ge \mu(M_{A,W})$ with equality iff $L$ is not special,
$s=r-1$, $r|g$,
$A$ is special and hence $A=L(-Z)$ for some effective divisor $Z$.
\qed\enddemo

\heading \S 2 Main Results \endheading

\demo{Proof of Main Theorem}
We need to construct a family parameterizing sequences:
$$0\to M_{L,V} \to V\otimes {\Cal O_C} \to L \to 0.$$
The line bundles $L$ are parameterized by the Poincare
bundle $\Cal P$ on the degree $d$ Jacobian  $J^d$.
Consider $X=\operatorname {Spec} (\operatorname {Sym} ({V^*}\otimes \Cal 
P) )$. There is a canonical section of $V^*\otimes \Cal P$ and this
gives rise to a sequence
$$0\to M_{\Cal P,V} \to V\otimes {\Cal O_C} \to \Cal P \to 0$$
over $X\times C$
which parameterizes the above sequences once we throw away
degenerate sequences (those where $V$ drops rank or does not span $L$).

We also need to parameterize sequences by using a `Poincare'
bundle on the space of $\rank =n$ stable bundles.  This is
a problem because there is no such bundle if $n$ and $d$ are
not coprime.  We can however, find finitely many Zariski
open subsets which cover the moduli space, and construct a Poincare
bundle over an etale cover of each cover so that the universal
map coincides with the projection map.  Using this we can 
construct our family as above.

Now let $\delta$ be the least integer such that 
$\rho (n,\delta,1)\ge 0$.
For $\delta \le d \le g+n$ we have a family of linebundles
$L$ which are spanned by global sections with $\dim \ho L)=n+1$.
And by Theorem 2 the dual span is generically stable.  
Now as above we get a family of sequences giving a spanning pair
and it's dual.  

But what if $d>g+n$?  Then we must use an
incomplete space of sections and Theorem 2 does not apply.
The way out of this mess is to start with the $\rank n$
vector bundles, show they are stable and spanned, and their
kernel will be a line bundle which must be stable.

Since $g\ge n^2-1$, $\rho (n,g+1,n)\ge 0$ and hence
$S_{\alpha}(n,g+1,n+1)$ is not empty.  Let $E$ be a
stable bundle with $(E,V)$ an element of that space.  
Consider a point $p\in C$. $E(p)\in S_{\alpha}(n,g+1+n,n+1)$.
So the latter space is non-empty.  It has a component on
which a general element is generically spanned and whose
spanned subbundle has no trivial summands.  
We now do a dimension count to show a general bundle is spanned.
If a generic bundle has spanned subbundle $F$ with 
$\deg(F) = (g+1+n - a)$, then the space of these bundles
has $\dim = \rho (n,g+1+n,n) - a(n+1)$. 
The dimension of bundles each subbundle can fit into
is $a(n)$ as found by counting possible elementary transformations.
So the dimension of bundles we started with is
$\rho (n,g+1+n,n) -a$.  Hence $a$ is zero because the
minimum  dimension is given by $\rho$ (assuming the space
is not empty).  In conclusion, the generic bundle is spanned.  
Furthermore, this technique works for all $d>g+n$.

Now that the family is constructed we use it to obtain
a birational map from $S_n = S_{\alpha}(n,d,n+1)$ to
$S_1 = S_{\alpha} (1,d,n+1)$.  There is a natural
map from $X$ to $S_1 \times S_n$; we call the image $S$.
Since the map is given by a dual span, if $x\in S$ and
$y \in S$ and the image of $x$ and $y$ correspond in
$S_1$ then they correspond in $S_n$, and vice versa.
So $S$ maps injectively into $S_1$ and $S_n$.  The
maps are then birational on the reduced schemes.

All that remains is the case $g=n=2$.  We do that
as an example.
\qed
\enddemo

\subheading{Example 1}
Let $g=n=2$.  If $\ho L)=3$, then $\deg (L) = 4=2g$.
By Theorem 2, $M_L$ is semistable but not stable.
In particular $M_{\omega_C} = \omega_C^*$ is a subbundle
of $M_L$ and hence $\omega_C$ is a quotient bundle of $M_L^*$.
Assume for the moment that $L \ne \omega_C^{\otimes 2}$,
and let $A=L\otimes \omega_C^*$.  There is a sequence
$$0 \to A \to M_L^* \to \omega_C \to 0.$$
$A$ is the only subbundle with $\mu(A)=\mu(M_L)$.
But $\ho A)=1$ and so $(\bho L)^*, M_L^*)$ is $\alpha$
stable.  So now let $A=\omega_C$.  We have a sequence
$$0 \to \omega_C \to M_L^* \to \omega_C \to 0.$$
If $M_L^*$ is spanned, so is the cokernel $\omega_C$.
That means the sequence is exact on global sections
and $\ho M_L^*)\ge 4$.  Furthermore, the endomorphisms
(given by scalar multiplication and surjection onto
the subbundle $\omega_C$ of $M_L^*$) has dimension $2$ if
the bundle is indecomposable, and $4$ otherwise.
In the first case we have a 3 dimensional family
of subspaces of the space of sections.  This is
acted upon by the group of endomorphism.  Modding
out by scalars we get a $2$ dimensional family
of stable pairs.  This is impossible because there
is only one dual span.  So the bundle decomposes.
Now there is a three dimensional family of subspaces
acted on by a 3 dimension group of automorphisms
and we get a unique pair (up to isomorphism).
Furthermore, the pair is spanned and hence any subbundle
has only a $1$ dimensional family of sections
(since the cokernel, which is spanned, has a two
dimensional space of sections).

Now the unspanned bundles all have subbundles $\omega_C$.
So it is easily seen that they are $\alpha$ unstable.

Now all of this gives us our isomorphism for $d=4$.
What about $d>4$.  We have no stable bundles to tensor
by a linebundle.  (And semistable does no good because
the determinantal loci is defined only for stable bundles.)
So we consider elementary transformations.
We make an ad hoc definition and call a primitive transformation
of $E$ at a reduced point p to be the kernel of a sequence:
$$0\to E_p \to E \to \Cal O_p \to 0.$$
We can dualize $E$ take a primitive transformation
and dualize again to get:
$$0 \to E \to E_p \to \Cal O_p \to 0.$$
By abuse of terminology and notation we call this
a primitive transformation.
If we can show a primitive transformation is stable,
then we can do the dimension count done in the proof
of the Main Theorem.

The proof that primitive transformations are stable
is given by Lange and Narasimhan \cite{4, Lemma 4.3}.
Each rank 2 vector bundle can be thought of as a
ruled surface.  The surface has a minimal section
(meaning a section with minimal self intersection $s$).
The self intersection is $>0$ (or $\ge 0$) iff
$E$ is stable (or semistable).  The minimal section
is not unique but with $g\ge 2$ there are only finitely
many unless $s\ge 2$ or $E$ is trivial.  

The point is that a primitive transformation raises
$s$ unless the transformation corresponds to a point
on a minimal section.  If $s=1$ or $0$, and $E$ is
not trivial, there are only finitely many minimal sections,
and therefor, $s$ rises for a general point 
and $E_p$ is stable.
If $s\ge 2$ it may drop, but only by one, which would
leave it stable.

The upshot is that primitive transformations show
$G_{\alpha}(2,d,3)$ is non-empty for $d\ge 5$.
And furthermore, the dimension count in the proof
of the Main Theorem applies to prove $S_{\alpha}(2,d,3)$
is also non-empty.

\heading \S 3 A Nonexistence Result \endheading

\proclaim{Theorem 3} Assume $C$ is general in the sense that if
$\rho (g,r,d)<0$, then $W^r_d = \emptyset$.
Then it follows that for
$\rho (g,n,n,d) <0$, then $W^n_{n,d} = \emptyset$.
\endproclaim

To prove this we need an elementary lemma.

\proclaim{Lemma 1} Let $E$ be a vector bundle of rank $n$.
$\rho (g,n,n,d) \ge 0$ is equivalent to the inequality 
$\mu(E)=\frac{d}{n} \ge 1+\frac{g}{n+1} 
\overset{\text{def} } \to = f(n)$.
\endproclaim

\demo{Proof of Lemma 1}
Use $\rho (g,n,n,d)=n^2(g-1)+1 - (n+1) (n+1-d+n(g-1))$
and solve for $\frac {d}{n}$.
\qed
\enddemo

\subheading {Remark 3} 
The function $f(n)$ is strictly decreasing.

\demo{Proof of Theorem 3}
Assume $\rho (g,n,n,d) <0$ and hence $\rho (g,n,1,d)<0$
which implies $W^n_{1,d}=\emptyset$.  
Now assume $E\in W^n_{n,d}$.
We have three cases:
\roster
\item $E$ spanned,
\item the global sections span a proper subbundle of the same rank, and
\item the global sections span a subbundle of smaller rank.
\endroster

Case 1.  If $E$ is spanned, then the dual span is a line bundle
$L$ with $L \in W^n_{1,d} = \emptyset$.  So this is impossible.

Case 2. Let $F$ be the proper subbundle spanned by the sections
of $E$.  We have two subcases.  A) $F$ has no trivial summands,
and B) $F=\oplus \Cal O_C \oplus G$ where $G$ has no trivial
summands.  Case A follows from the proof of Case 1.  As for 
$B$ we note that $h^0(C,G)\ge \operatorname{rank}(G)+1$.
But stability of $E$ and the fact that $f(n)$ is strictly
decreasing  gives
$$\mu(G) \le \mu(E) < f(n) < f(\operatorname{rank}(G) ).$$
So by induction, case 2B follows.

Case 3.  Argue as in 2B.
\qed
\enddemo


\enddocument